\documentclass[11pt]{article}

\usepackage{latexsym,epsfig,amssymb,amsmath,amsthm,color,url,bbm,pdfsync}
\usepackage{graphicx,float,bm}
\usepackage{subcaption}
\usepackage{dsfont}

\numberwithin{equation}{section}
\allowdisplaybreaks
\usepackage[round,comma]{natbib}

\newtheorem{lemma}{Lemma}[section]

\newtheorem{proposition}[lemma]{Proposition}
\newtheorem{remark}[lemma]{Remark}
\newtheorem{corollary}[lemma]{Corollary}
\newtheorem{definition}[lemma]{Definition}

\title{Characterizing extremal dependence on a hyperplane}
\author{
Phyllis Wan\footnote{Erasmus University Rotterdam; Econometric Institute, Burg.\ Oudlaan 50, 3062 PA Rotterdam, the Netherlands;
email: wan@ese.eur.nl}
}
\date{}

\begin{document}
\maketitle

\begin{abstract}

In this paper, we characterize the extremal dependence of $d$ asymptotically dependent variables by a class of random vectors on the $(d-1)$-dimensional hyperplane perpendicular to the diagonal vector $\mathbf1=(1,\ldots,1)$.  This translates analyses of multivariate extremes to that on a linear vector space, opening up possibilities for applying existing statistical techniques that are based on linear operations.  As an example, we demonstrate obtaining lower-dimensional approximations of the tail dependence through principal component analysis.  Additionally, we show that the widely used H\"usler-Reiss family is characterized by a Gaussian family residing on the hyperplane.

\end{abstract}


\section{Introduction}

%
%
%
Let $\mathbf{Y} = (Y_1,\ldots,Y_d)$ be a random vector with continuous marginal cdf's $F_1,\ldots,F_d$.  There are two common approaches in the literature to geometrically characterize its tail dependence.
\begin{enumerate}
\item
	{\bf Angular component}: Transform the marginals to standard Pareto with $\tilde{\mathbf{X}} = (\tilde{X}_1,\ldots,\tilde{X}_d) = \left(1/\{1-F_1(Y_1)\},\ldots,1/\{1-F_d(Y_d)\}\right)$.  Conditioning on the norm of $\tilde{\mathbf{X}}$ being large for a pre-specified norm $\|\cdot\|$,
	$$
		\left.\frac{\tilde{\mathbf{X}}}{r} \ 	\right|\ \{\|\tilde{\mathbf{X}}\| >r\} \overset{d}\to R\cdot \Theta, \quad  r\to\infty.
	$$
	Here $\Theta$, the {\it angular component}, resides on the positive unit sphere $\{\mathbf{v}\in [0,\infty)^d| \|\mathbf{v}\|=1\}$ and $R$ is a standard Pareto variable independent of $\Theta$. This follows from the framework of {\it multivariate regular variation}, see e.g., \cite{R2007}.
\item
	{\bf Spectral random vector}:  Transform the marginals to standard exponential with $\mathbf{X} = (X_1,\ldots,X_d) = \left(-\log\{1-F_1(Y_1)\},\ldots,-\log\{1-F_d(Y_d)\}\right)$. Conditioning on the maximum component of $\mathbf{X}$ being large, 
	$$
		\mathbf{X} - r\cdot \mathbf{1} \  \left| \  \{\max(\mathbf{X}) > r\} \right. \overset{d}\to E\cdot \mathbf1 + \mathbf{S}, \quad  r\to\infty.
	$$
	Here $\mathbf{S}$, the {\it spectral random vector}, resides on the space $\{\mathbf{v}\in [-\infty,0]^d|\max(\mathbf{v}) = 0\}$ and $E$ is a unit exponential ($\text{Exp}(1)$) variable independent of  $\mathbf{S}$.  This follows from the framework of {\it multivariate peak-over-threshold}, see \cite{RT2006} and \cite{RSW2018}.
\end{enumerate}

The two characterizations are connected as the latter is equivalent to the former with the $L_\infty$-norm.  While either $\Theta$ or $\mathbf{S}$ serves to summarize the extremal dependence of $\mathbf{X}$, both possess nonlinear supports that induce intrinsic dependency between dimensions.  This poses nontrivial constraints for the construction of statistical models and their inference when it comes to studying the tails.

In this paper, we focus on a random vector $\mathbf{X}$ with $\text{Exp}(1)$-like marginals and instead condition on the component mean $\bar{X} = d^{-1}\sum_{k=1}^d X_k$ being large.  We consider the scenario where $\mathbf{X}$ has asymptotically dependent components such that they are simultaneously large in the tail.  This translates to each component of the spectral random vector $\mathbf{S}$ having no mass on $-\infty$, see Section~\ref{sec:background} for details.  We show that
$$
	\mathbf{X} - r\cdot \mathbf{1} \  \left|\  \{\bar{X} > r\} \right. \overset{d} \to E\cdot \mathbf1 + \mathbf{U}, \quad \text{as } r\to\infty,
$$
where $\mathbf{U}$ belongs to the class of random vectors $\mathcal{U} = \{\mathbf{U} \in \mathbf1^\perp|E\{e^{\max(\mathbf{U})}\}<\boldsymbol\infty\}$, with $\mathbf1^\perp:= \{\mathbf{v}|\mathbf{v}^T\mathbf1=0\}$ being the hyperplane perpendicular to the diagonal vector $\mathbf1$, and $E$ is an $\text{Exp}(1)$ variable independent of $\mathbf{U}$.  We term $\mathbf{U}$ the {\it profile random vector} and propose it as an alternative summary for extremal dependence.   

We point out two attractive properties in this proposal.  First, the class $\mathcal{U}$ resides on a linear vector space and is closed under finite addition and scalar multiplication.  This allows existing statistical techniques based on linear operations to be straightforwardly adapted for extremes.  As an example, we illustrate the application of principal component analysis to achieve lower-dimensional approximations of the tails.
Second, the H\"usler-Reiss family \citep{HR1989}, the class of nontrivial tail dependence of Gaussian triangular arrays, is characterized by Gaussian profile random vectors.  Despite being one of the most widely used parametric models for extremes, the analytical form of H\"usler-Reiss models is not easy to handle mathematically.  Using profile random vectors, analyses of H\"usler-Reiss models can be translated to that of Gaussian models on the hyperplane $\mathbf1^\perp$.

\subsection*{Notation}

Boldface symbols are used to denote vectors, usually of length $d$.  We write $\mathbf0=(0,\ldots,0)$ and $\mathbf1=(1,\ldots,1)$, where the lengths of the vector may depend on the context.  The maximum component and component mean of $\mathbf{x}=(x_1,\ldots,x_d)$ are denoted by  $\max(\mathbf{x}) = \max(x_1,\ldots,x_d)$ and $\bar{x} = d^{-1}\sum_{k=1}^d x_k$, respectively.  When applied to vectors, mathematical operations, such as addition, multiplication, exponentiation, maximum and minimum, are taken to be component-wise.  Lastly, $\mathbf{1} ^\perp:= \{\mathbf{v}|\mathbf{v}^T\mathbf1=0\}$ is used to denote the hyperplane perpendicular to the vector $\mathbf1$.


\section{Background on multivariate extremes} \label{sec:background}

    Let $\mathbf{X}$ be a random vector in $\mathbb{R}^d$.  To study the tail of $\mathbf{X}$, a common assumption is that there exist sequences of normalizing vectors $\{\mathbf{a}_n\}$ and $\{\mathbf{b}_n\}$, such that for $\mathbf{X}_1,\ldots,\mathbf{X}_d \overset{iid} \sim \mathbf{X}$,
\begin{equation}\label{eq:bm}
	\lim_{n\to\infty} \textup{pr}\left( \frac{\max_{i=1,\ldots,n}\mathbf{X}_i - \mathbf{b}_n}{\mathbf{a}_n} \le \mathbf{x}\right) = G(\mathbf{x}).
\end{equation}
Such a limit distribution $G$ with non-degenerate marginals is referred to as a {\it generalized extreme value distribution} and $\mathbf{X}$ is said to be in the domain of attraction of $G$.  The marginals of $G$  can be parametrized by
$$
	G_k(x_k) = \exp\left[ - \left\{1 + \gamma_k (x_k - \mu_k)/\alpha_k\right\}^{-\frac{1}{\gamma_k}}\right], \quad 1 + \gamma_k (x_k - \mu_k)/\alpha_k>0,
$$
where $\gamma_k,\mu_k \in \mathbb{R}$ and $\alpha_k>0$.  In the case where $\gamma_k=0$, $G_k(x_k)$ is interpreted as the limit $G_k(x_k) = \exp[-\exp\{-(x_k-\mu_k)/\alpha_k\}]$.  The dependence structure of $G$, on the other hand, cannot be parametrized and may be complex.  For background on multivariate generalized extreme value distributions and their domains of attraction, see e.g.,~\cite{dHF2006}.

To focus exclusively on extremal dependence, we assume that the marginals of $\mathbf{X}$ are transformed to similar scales defined as follows.
\begin{definition} \label{def:exp}
Define $\mathcal{X}$ to be the class of random vectors $\mathbf{X}$ such that for $\mathbf{X}_1,\ldots,\mathbf{X}_d \overset{iid} \sim \mathbf{X}$,
\begin{equation}\label{eq:stbm}
	\lim_{n\to\infty} \textup{pr}\left\{\max_{i=1,\ldots,n}\mathbf{X}_i - \log(n) \cdot \mathbf1 \le \mathbf{x}\right\} = G(\mathbf{x}),
\end{equation}
where the marginals of $G$ follow the Gumbel distribution $G_k(x_k) = \exp[-\exp\{-(x_k-\mu_k)\}]$.  
\end{definition}
The class $\mathcal{X}$ describes random vectors with marginals of similar scales as $\text{Exp}(1)$.  An arbitrary random vector can be transformed to an element of $\mathcal{X}$ by standardizing its marginals to $\text{Exp}(1)$.

Following elementary calculation from \eqref{eq:stbm} \citep{RT2006}, 
\begin{equation} \label{eq:mpot:std}
	\mathbf{X} - r \cdot \mathbf1\  \left|\  \{\max(\mathbf{X})\ge r\} \right. \overset{d}\to \mathbf{Z}, \quad r \to\infty,
\end{equation}
where $\mathbf{Z}$ has distribution function 
\begin{equation} \label{eq:mgpd:df}
	\textup{pr}(\mathbf{Z} \le \mathbf{z}) = \frac{\ln G(\mathbf{z} \wedge \mathbf0) - \ln G(\mathbf{z})}{\ln G(\mathbf0)},
\end{equation}
and is referred to as a {\it standardized multivariate generalized Pareto distribution} \citep{RSW2018}.
\begin{definition} \label{def:smgpd}
A random vector $\mathbf{Z}$ is a {\it standardized multivariate generalized Pareto distribution} if there exists $\mathbf{X} \in \mathcal{X}$ such that \eqref{eq:mpot:std} holds.
\end{definition}

In this paper, we focus on the scenario where the components of $\mathbf{X}$ are {\it asymptotically dependent}.
\begin{definition} \label{def:smgpd}
A random vector $\mathbf{X} \in \mathcal{X}$ or its corresponding standardized multivariate generalized Pareto distribution $\mathbf{Z}$ is said to have asymptotically dependent components if $\textup{pr}(Z_k > -\infty)=1$, $k=1,\ldots,d$.
\end{definition}
Intuitively, asymptotically dependent components are large simultaneously in the tail.  This scenario serves as a foundation for studying more complicated extremal dependence.  A generic tail dependence structure can be constructed via a mixture model, where each factor is asymptotically dependent on a selection of components and degenerate on the rest, such that each factor can be modelled by a lower-dimensional multivariate generalized Pareto distribution with asymptotically dependent components.  Such a framework was proposed for multivariate generalized Pareto distributions in \cite{MKS2024} and we refer the readers to the reference therein for earlier work on the detection and modelling of asymptotic dependence in general.

The following proposition, adapted from Theorem~7 of \cite{RSW2018}, shows that $\mathbf{Z}$ can be represented by a random vector on the L-shaped space $\{\mathbf{v}| \max(\mathbf{v})=0\}$. 

\begin{proposition} \label{prop:rep:mgp}

Let $\mathcal{S}$ be the class of random vectors $\mathbf{S} \in (-\infty,0]^d$ such that 
$\textup{pr}\{\max(\mathbf{S})=0\}=1$ and
$\textup{pr}(S_k > -\infty)=1$, $k=1,\ldots,d$.  Let $\mathbf{Z}$ be a standardized multivariate generalized Pareto distribution with asymptotically dependent components.  
Then 
\begin{equation} \label{eq:mgp:s}
	\mathbf{Z} \overset{d}{=} E \cdot \mathbf1 + \mathbf{S},
\end{equation}
where $\mathbf{S} \in \mathcal{S}$ and $E$ is an $\text{Exp}(1)$ variable independent of $\mathbf{S}$.  Conversely, any $\mathbf{S} \in \mathcal{S}$ characterizes a standardized multivariate generalized Pareto distribution $\mathbf{Z}$ with asymptotically dependent components through \eqref{eq:mgp:s}. 
\end{proposition}

The random vector $\mathbf{S}$ is referred to as the {\it spectral random vector} associated with $\mathbf{Z}$. An illustration of the domains of $\mathbf{Z}$ and $\mathbf{S}$ is shown in Figure~\ref{fig:1}(a).  Effectively, the spectral random vector is the limit 
$$
	\mathbf{X} - \max(\mathbf{X}) \cdot \mathbf1\  \left|\  \{\max(\mathbf{X})\ge r\} \right. \overset{d}\to \mathbf{S}, \quad r\to\infty,
$$
representing the tail of $\mathbf{X}$ being diagonally projected onto the L-shaped space $\{\mathbf{v}|\max(\mathbf{v})=0\}$.


\section{Diagonal peak-over-threshold and profile random vectors} \label{sec:dmgp}

\subsection{Diagonal peak-over-threshold}

In this section, we consider a different peak-over-threshold framework.   We propose thresholding the tail based on the component mean instead of the maximum component.

\begin{proposition} \label{prop:dpot}
Given $\mathbf{X} \in \mathcal{X}$, let $\mathbf{Z}$ be its corresponding standardized multivariate generalized Pareto distribution.  Assume that $\mathbf{X}$ and $\mathbf{Z}$ have asymptotic dependent components.  Then
$$
	\mathbf{X} - r \cdot \mathbf1\ \left|\   \bar{X} \ge r \right. \overset{d}\to \mathbf{Z}^*, \quad r \to\infty,
$$
where
\begin{equation} \label{eq:zstar}
	\mathbf{Z}^*: \overset{d} = \mathbf{Z}\ \left|\ \{\mathbf{Z}^T\mathbf1 \ge0\} \right. .
\end{equation}
\end{proposition}

We refer to the limit distribution $\mathbf{Z}^*$ as a {\it diagonal multivariate generalized Pareto distribution}.  If a pair of standardized and diagonal multivariate generalized Pareto distributions $(\mathbf{Z},\mathbf{Z}^*)$ satisfies \eqref{eq:zstar}, then we say they are associated.
\begin{remark}
This paper does not focus on the scenario when a random vector has asymptotically independent components, where the components of $\mathbf{Z}$ have mass on $-\infty$ and $\{\mathbf{Z}^T\mathbf1>0\}$ may have probability 0.  Here we present a small illustration of what could happen.  Consider $\mathbf{X} = (X_1,X_2)$ with $\text{Exp}(1)$ margins and $\mathbf{Y} = (Y_1,Y_2) = (e^{X_1},e^{X_2})$ with standard Pareto margins.  Projecting the tail of $\mathbf{X}$ onto $\mathbf1^\perp=\{(x_1,x_2)|x_1+x_2=0\}$ is equivalent to projecting the tail of $\mathbf{Y}$ onto $\{(y_1,y_2)|(y_1y_2)^{1/2} =1\}$.  In the case where $X_1$ and $X_2$ (hence $Y_1$ and $Y_2$) are asymptotically independent, the projection reveals the hidden regular variation between $Y_1$ and $Y_2$  \citep{MR2004}. 
\end{remark}


%


\subsection{Profile random vectors}

\begin{figure}[t]
	\centering
	\begin{subfigure}{0.24\textwidth}
  	\includegraphics[width=\textwidth]{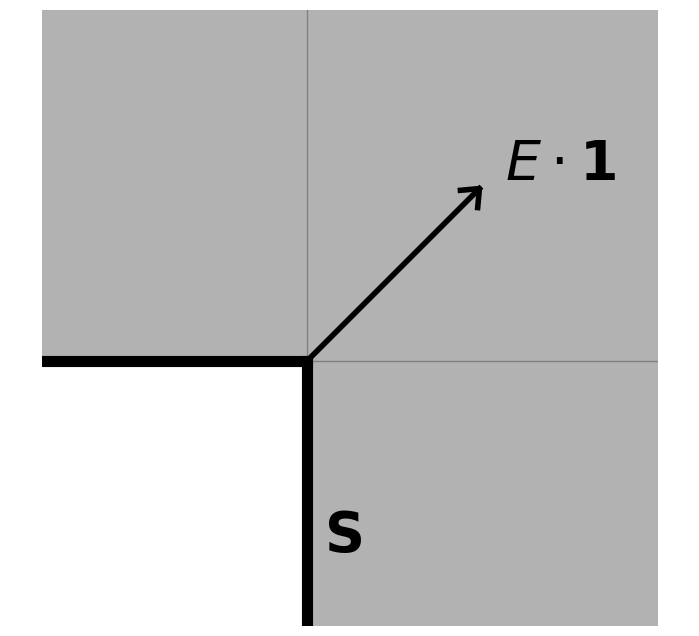}
	\caption{}
	\end{subfigure}
	\begin{subfigure}{0.24\textwidth}
  	\includegraphics[width=\textwidth]{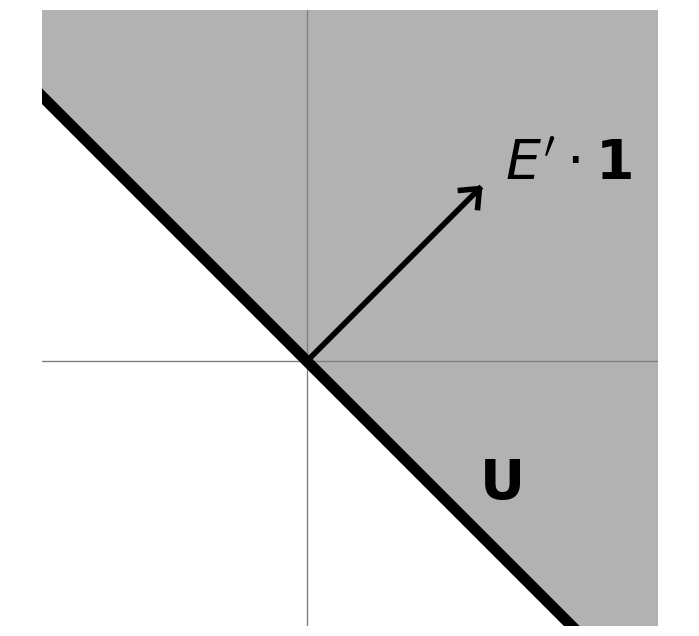}
	\caption{}
	\end{subfigure}
	\begin{subfigure}{0.24\textwidth}
  	\includegraphics[width=\textwidth]{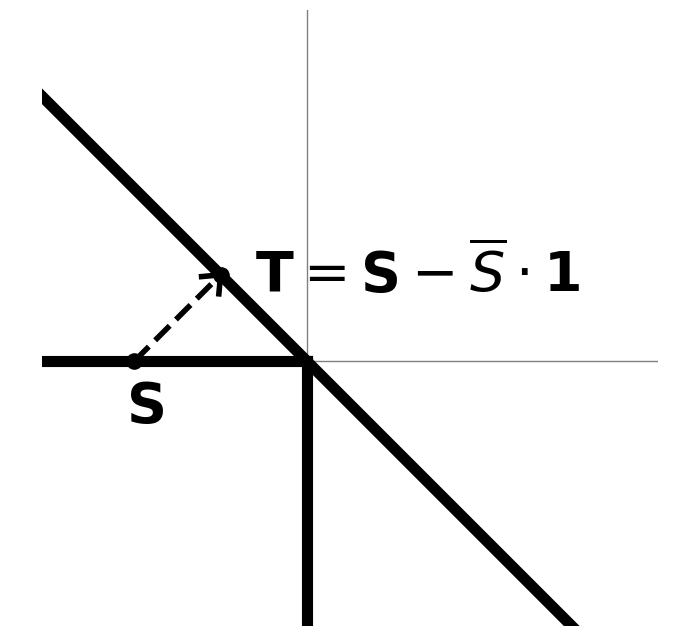}
	\caption{}
	\end{subfigure}
	\begin{subfigure}{0.24\textwidth}
  	\includegraphics[width=\textwidth]{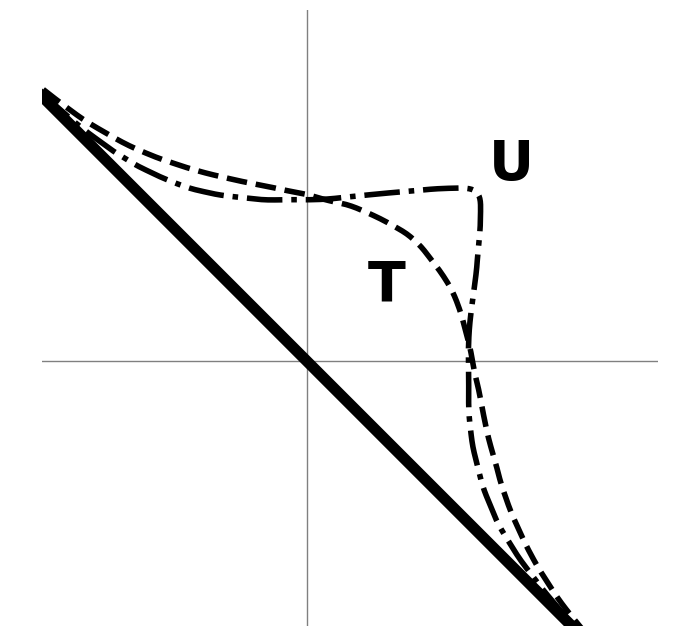}
	\caption{}
	\end{subfigure}
	\caption{Plot (a): Domain of $\mathbf{S}$ (thick line), direction of $E \cdot \mathbf1$ (arrow), and domain of $\mathbf{Z} \overset{d}= \mathbf{S} + E \cdot \mathbf1$ (shaded area); Plot (b): Domain of $\mathbf{U}$ (thick line), direction of $E' \cdot \mathbf1$ (arrow), and domain of $\mathbf{Z}^* \overset{d}= \mathbf{U} + E' \cdot \mathbf1$ (shaded area); Plot (c): Example of projection from $\mathbf{S}$ to $\mathbf{T} = \mathbf{S} - \bar{S} \cdot \mathbf1$;  Plot (d): densities of a pair of associated $\mathbf{T}$ and $\mathbf{U}$ on $\mathbf1^\perp$.}
	\label{fig:1}
\end{figure}

We now establish the links between associated standardized and diagonal multivariate generalized Pareto distributions $\mathbf{Z}$ and $\mathbf{Z}^*$.  Instead of the spectral random vector $\mathbf{S}$, $\mathbf{Z}$ can be equivalently characterized by the projection of $\mathbf{S}$ onto the hyperplane $\mathbf{1}^{\perp}$,
\begin{equation} \label{eq:t:generator}
	\mathbf{T} := \mathbf{S} - \bar{S} \cdot \mathbf1 \in \mathbf1^\perp,
\end{equation}
since $\mathbf{S}$ can be retrieved from $\mathbf{T}$ via $\mathbf{S} = \mathbf{T} - \max(\mathbf{T}) \cdot \mathbf1$.  The projection from $\mathbf{S}$ to $\mathbf{T}$ is illustrated in Figure~\ref{fig:1}(c).  The following proposition characterizes the distribution of $\mathbf{Z}^*$ through $\mathbf{T}$.

\begin{proposition} \label{prop:u}
Let $\mathbf{Z}$ and $\mathbf{Z}^*$ be associated standardized and diagonal multivariate generalized Pareto distributions.  Let $\mathbf{S}$ be the spectral random vector of $\mathbf{Z}$ and let $\mathbf{T}$ be as defined in \eqref{eq:t:generator}.  Then $\mathbf{Z}^*$ has stochastic representation
$$
	\mathbf{Z^*} \overset{d}= E' \cdot \mathbf1 + \mathbf{U},
$$
where $\mathbf{U} \in \mathbf1^{\perp}$ and $E'$ is an $\text{Exp}(1)$ variable is independent of $\mathbf{U}$.  The distribution of $\mathbf{U}$ is given by 
$$
	\mathbf{U} :\overset{d}= \mathbf{T}\ \left| \ \{\max(\mathbf{T}) \le E\} \right.,
$$
where $E$ is an $\text{Exp}(1)$ variable is independent of $\mathbf{T}$.
\end{proposition}
We refer to  $\mathbf{U}$ as the {\it profile random vector} of $\mathbf{Z}^*$. An illustration of the domains of $\mathbf{Z}^*$ and $\mathbf{U}$ is shown in Figure~\ref{fig:1}(b).   We say that a pair of spectral and profile random vectors $\mathbf{S}$ and $\mathbf{U}$ are associated if the corresponding $\mathbf{Z}$ and $\mathbf{Z}^*$ are associated.  

\begin{remark}
The notation of $\mathbf{T}$ and $\mathbf{U}$ is inherited from \cite{RSW2018}, which proposed that a spectral random vector $\mathbf{S}$ could be generated from a random vector $\mathbf{T}$ via $\mathbf{S} = \mathbf{T} - \max(\mathbf{T}) \cdot \mathbf1$, or from a random vector $\mathbf{U}$ via $\textup{pr}(\mathbf{S} \in \cdot) = E\left[\mathds{1}_{\{\mathbf{U} - \max(\mathbf{U}) \in \cdot \}}e^{\max(\mathbf{U})}\right]/E\left\{e^{\max(\mathbf{U})}\right\}$. It follows that our $\mathbf{T}$ and $\mathbf{U}$ corresponds to the unique such random vectors on $\mathbf1^\perp$. 
\end{remark}

\begin{corollary}  \label{cor:moments}
Let $\mathbf{S}$ and $\mathbf{U}$ be associated spectral and profile random vectors.  Let $\mathbf{T}$ be as defined in \eqref{eq:t:generator}.  Then $E\{e^{\max(\mathbf{U})}\} < \infty$.
\end{corollary}

We link the distributions of $\mathbf{T}$ and $\mathbf{U}$ by conditioning on their respective maximum components.
\begin{proposition} \label{prop:t:u}
Let $\mathbf{S}$ and $\mathbf{U}$ be associated spectral and profile random vectors.  Let $\mathbf{T}$ be as defined in \eqref{eq:t:generator}.  Then 
$$
	 \mathbf{U} \ \left|\ \{\max(\mathbf{U})=s\} \right. \overset{d}= \mathbf{T} \ \left|\ \{\max(\mathbf{T})=s\} \right., \quad  s \ge 0.
$$
Given $\max(\mathbf{T})$, the distribution of $\max(\mathbf{U})$ can be obtained from
\begin{equation} \label{eq:maxu}
	\textup{pr}\{\max(\mathbf{U}) \le s\} = \frac{\int_0^s \textup{pr}\{\max(\mathbf{T}) \le t\}  e^{-t} dt + e^{-s} \textup{pr}\{\max(\mathbf{T}) \le s\}}{E\left\{e^{-\max(\mathbf{T})}\right\}}, \quad s \ge 0.
\end{equation}
Conversely, given $\max(\mathbf{U})$, the distribution of $\max(\mathbf{T})$ can be obtained from
\begin{equation} \label{eq:maxt}
	\textup{pr}\{\max(\mathbf{T}) \le s\} = \frac{e^{s} \textup{pr}\{\max(\mathbf{U}) \le s\} - \int_0^s \textup{pr}\{\max(\mathbf{U}) \le  t \}  e^{t} dt}{E\left\{e^{\max(\mathbf{U})}\right\}}, \quad s \ge 0.
\end{equation}
\end{proposition}

In the case where $\max(\mathbf{T})$ and $\max(\mathbf{U})$ are absolutely continuous, the link is simplified via density functions.

\begin{corollary}  \label{cor:t:u:2}
Let $\mathbf{S}$ and $\mathbf{U}$ be associated spectral and profile random vectors.  Let $\mathbf{T}$ be as defined in \eqref{eq:t:generator}.   If $\max(\mathbf{T})$ is absolutely continuous and admits density $f_{\max(\mathbf{T})}$, then $\max(\mathbf{U})$ is absolutely continuous with density
$$
	f_{\max(\mathbf{U})}(s) = \frac{1}{E\left\{e^{-\max(\mathbf{T})}\right\} } \  f_{\max(\mathbf{T})}(s)\  e^{-s}.
$$
Conversely, if $\max(\mathbf{U})$ is absolutely continuous and admits density $f_{\max(\mathbf{U})}$, then $\max(\mathbf{T})$ is absolutely continuous with density
$$
	f_{\max(\mathbf{T})}(s) =  \frac{1}{E\left\{e^{\max(\mathbf{U})}\right\} } \ f_{\max(\mathbf{U})}(s)\ e^{s}.
$$
\end{corollary}

An illustration of the densities of a pair of associated $\mathbf{T}$ and $\mathbf{U}$ is shown in Figure~\ref{fig:1}(d).  Effectively, conditioning $\mathbf{T}$ on $\max(\mathbf{T}) \le E$ shrinks its `radius' $\max(\mathbf{T})$ such that $\mathbf{U}$ is more concentrated around $\mathbf0$.

On the other hand, $\mathbf{T}$ can be obtained from $\mathbf{U}$ through the following stochastic limit.

\begin{corollary} \label{cor:t:u:1}
Let $\mathbf{S}$ and $\mathbf{U}$ be associated spectral and profile random vectors.  Let $\mathbf{T}$ be as defined in \eqref{eq:t:generator}.  
Given an $\text{Exp}(1)$ variable $E'$ independent of $\mathbf{U}$,
$$
	\mathbf{U}\ \left|\ \{\max(\mathbf{U}) \ge r-E'\} \right. \overset{d} \to \mathbf{T}, \quad r\to\infty.
$$

\end{corollary}

Finally, let $\mathcal{U}$ be the class of random vector $\mathbf{U} \in \mathbf1^\perp$ such that $E\{e^{\max(\mathbf{U})}\}<\infty$. The following proposition shows that any $\mathbf{U} \in \mathcal{U}$ is a profile random vector.
\begin{proposition} \label{prop:rv}

Given any $\mathbf{U} \in \mathcal{U} = \{\mathbf{U} \in \mathbf1^\perp| E\{e^{\max(\mathbf{U})}\}<\infty\}$,  let $E$ be an $\text{Exp}(1)$ variable independent of $\mathbf{U}$.  Then the random vector defined by
$$
	\mathbf{X} :\overset{d}= E \cdot \mathbf1 + \mathbf{U}.
$$
satisfies $\mathbf{X} \in \mathcal{X}$.  In particular, the profile random vector associated with its tail is $\mathbf{U}$.
\end{proposition}


\section{Principal component analysis} \label{sec:pca}

The class of profile random vectors $\mathcal{U}$ resides on a linear vector space and is closed under finite addition and scalar multiplication.  This provides a context to apply statistical analysis based on linear techniques for extremes.  In the following, we illustrate the application of principal component analysis as an example.

Principal component analysis is a classical technique for finding lower-dimensional representations of a random vector while retaining most of its variability.  Given a centered random vector $\mathbf{X} \in \mathbb{R}^d$, it identifies the linear subspace $\mathcal{S}_p^*\subset \mathbb{R}^d$ of dimension $p < d$ such that the $L_2$-distance between $\mathbf{X}$ and its projection onto $\mathcal{S}_p^*$, $\Pi_{\mathcal{S}_p^*}\mathbf{X}$, is minimized.
This is achieved by computing the orthonormal eigenvectors $\mathbf{v}_1,\ldots,\mathbf{v}_d$ of the covariance matrix $E(\mathbf{X}\mathbf{X}^T)$ with ordered eigenvalues $\lambda_1\ge \ldots\ge \lambda_d\ge0$.  The optimal subspace $\mathcal{S}_p^*$ is the span of $\mathbf{v}_1,\ldots,\mathbf{v}_p$ and the best $p$-dimensional approximation of $\mathbf{X}$ is its projection onto the first $p$ principal components $\Pi_{\mathcal{S}_p^*}\mathbf{X} = \Pi_{\mathbf{v}_1}\mathbf{X} + \cdots + \Pi_{\mathbf{v}_p}\mathbf{X}$.

Previous attempts to apply principal component analysis to extremes focused on the angular component $\Theta$, see \cite{CT2019} and \cite{DS2021}.  However, since $\Theta$ resides on the unit sphere, any lower-dimensional approximations of $\Theta$ via principal component analysis no longer result in angular components.

We now point out that instead it is natural to apply principal component analysis to the profile random vector $\mathbf{U}$.  Without loss of generality, assume that $\mathbf{U}$ is centered.  Otherwise subtract by its mean.  Given $E[e^{\max(\mathbf{U})}]<\infty$, the covariance matrix $E(\mathbf{U}\mathbf{U}^T)$ always exists.  Since $\mathbf{U} \in \mathbf1^\perp$, the last eigenvector $\mathbf{v}_d$ is proportional to $\mathbf1$ with eigenvalue $\lambda_d=0$.
Every other eigenvector $\mathbf{v}_k$ is perpendicular to $\mathbf{v}_d$ and hence satisfies $\mathbf{v}_k \in \mathbf1^\perp$. For any $p < d-1$, $\Pi_{\mathcal{S}_p^*}\mathbf{U} = \Pi_{\mathbf{v}_1}\mathbf{U} + \cdots + \Pi_{\mathbf{v}_p}\mathbf{U} \in \mathbf1^\perp$ is a lower-dimensional approximation of $\mathbf{U}$, defines a profile random vector and induces an extremal dependence structure with lower dimensions.

%

In the conventional principal component analysis, the discarded principal components describe directions along which the variation of the data is minimized.  In the context of profile random vectors, the discarded principal components describe directions where extremal dependence is strong enough to be approximated by complete dependence.  To see this, consider the trivial case where $\mathbf{U}$ can be approximated by $\mathbf{0}$. Then the diagonal multivariate generalized Pareto distribution $\mathbf{Z}^* \overset{d}= E\cdot \mathbf1$ lies on the vector $\mathbf1$, meaning that all components are completely dependent in the tail.


\section{Gaussian profile random vectors} \label{sec:models}

Any parametric family on $\mathcal{U}$ induces a parametric family for profile random vectors. 
In particular, a class of Gaussian profile random vectors results in the family of H\"usler-Reiss models, the class of distributions describing the non-trivial tail limit of Gaussian triangular arrays \citep{HR1989}.  On this we elaborate in the following.

Any pair of components of a Gaussian random vector are asymptotically independent unless being perfectly collinear \citep{S1960}.  To construct nontrivial tail dependence, consider instead a Gaussian triangular array $\mathbf{X}_{i}^{(n)} \sim N({\bf0}, \Sigma^{(n)})$, $i=1,\ldots,n$, where $\Sigma^{(n)}_{kk}=1$, $k=1,\ldots,d$, and the off-diagonal elements of $\Sigma^{(n)}$ converge to 1 such that
$$
	\log(n) \cdot (\mathbf1\mathbf1^T - \Sigma^{(n)}) \to \Gamma = \left(\Gamma_{ij}\right)_{1\le i,j\le d}.
$$
The matrix $\Gamma$ satisfies $\Gamma_{ij} = E(W_i-W_j)^2$ for some centered multivariate Gaussian random vector $\mathbf{W}=(W_1,\ldots,W_d)$ and is called the {\it variogram} of $\mathbf{W}$.  Let $\phi(\cdot)$ be the density function of a standard normal variable and let $b_n$ be the solution to the equation $b_n=\phi(b_n)$.  Then
$$
	\lim_{n\to\infty} \textup{pr}\left\{b_n \cdot \left(\max_{i=1,\ldots,n}\mathbf{X}^{(n)}_i - b_n \cdot \mathbf1\right)\le \mathbf{x}\right\} = G_\Gamma(\mathbf{x}).
$$
The limiting distribution $G_\Gamma$ has standard Gumbel marginals and the generalized multivariate extreme value distribution associated with the {\it H\"usler-Reiss model parametrized by $\Gamma$}.  The standard multivariate generalized Pareto distribution $\mathbf{Z}$ for the H\"usler-Reiss model is defined accordingly from $G_\Gamma$ via \eqref{eq:mgpd:df}.

The following result shows that the profile random vector of a H\"usler-Reiss model is Gaussian.

\begin{proposition} \label{prop:hr}

The profile random vector of the H\"usler-Reiss model parametrized by $\Gamma$ is 
$$
	\mathbf{U} \sim N\left(\boldsymbol{\mu},\Sigma\right),
$$
where
$$
	\Sigma := -\frac{1}{2}\left(I - \frac{\mathbf{1}\mathbf{1}^T}{d}\right) \Gamma \left(I - \frac{\mathbf{1}\mathbf{1}^T}{d}\right), 
$$
and 
\begin{equation} \label{eq:hr:mu}
	\boldsymbol\mu:= -\frac{1}{2} \left\{\textup{diag}(\Sigma) - \overline{\textup{diag}(\Sigma)} \cdot \mathbf1\right\}.
\end{equation}
\end{proposition}
Proposition~\ref{prop:hr} was independently derived in an unpublished manuscript by Johan Segers in 2019.  In the special case where $\Sigma$ is of rank $(d-1)$, this result can be seen from Proposition~3.6 of \cite{HES2024}. In recent literature on H\"usler-Reiss models, $\Gamma$ is often assumed to be the variogram of a full-rank Gaussian vector such that the resulting multivariate generalized Pareto distribution $\mathbf{Z}$ admits a density.  The resulting $\Sigma$ is then of rank $(d-1)$ and its pseudo-inverse embeds information on the conditional independence  in the tail \citep{HES2024,WZ2023}.   In contrast, the result in Proposition~\ref{prop:hr} applies to H\"usler-Reiss models of all ranks.

Proposition~\ref{prop:hr} can facilitate the use of H\"usler-Reiss models in multiple aspects.  First, random vectors with H\"usler-Reiss extremal dependence can be generated directly via Proposition~\ref{prop:rv}, instead of being approximated via a triangular array by definition.  Second, lower-dimensional approximation of a H\"usler-Reiss model can be obtained by approximating $\mathbf{U}$ with a lower-dimensional Gaussian vector through principal component analysis.  Lastly, inference on H\"usler-Reiss parameters may be carried out through likelihood methods by diagonally thresholding the tail observations.

Finally, $\boldsymbol\mu$ and $\Sigma$ in the H\"usler-Reiss profile random vector are linked in \eqref{eq:hr:mu}.  This comes from the assumption that the $\mathbf{X}_i^{(n)}$'s in the triangular array have identical margins.  By relaxing this assumption, the family of H\"usler-Reiss tails can be extended by considering the class of all Gaussian profile random vector on $\mathbf1^\perp$.

\section*{Acknowledgement}

The author thank the associate editor and two anonymous reviewers for helpful comments, as well as Anja Jan\ss en, Chen Zhou, and other participants of Oberwolfach Workshop {\it Mathematics, Statistics, and Geometry of Extreme Events in High Dimensions} (2024) for extensive discussions.  The illustrations in Figure 1 were provided by Jasper Velthoen.  This research was supported by a Veni grant from the Dutch Research Council.


\appendix

\section{Proofs} \label{sec:proofs}

This section contains proofs for the propositions and the corollaries in the paper.  All equation numbers refer to those in the manuscript.

\begin{proof}[Proof of Proposition~\ref{prop:rep:mgp}]

From Theorem~7 of \cite{RSW2018}, $\mathbf{Z}$ admits the stochastic representation
$$
	\mathbf{Z} \overset{d}{=} E \cdot \mathbf1 + \mathbf{S},
$$
where $\mathbf{S}$ satisfies $\textup{pr}(\max(\mathbf{S})=0)=1$ and $E$ is an $\text{Exp}(1)$ variable independent of $\mathbf{S}$.   Since $\mathbf{Z}$ has asymptotically dependent components, $\textup{pr}(Z_k > -\infty)=1$, $k=1,\ldots,d$, and hence $\textup{pr}(S_k > -\infty)=1$, $k=1,\ldots,d$.  Consequently, $\mathbf{S} \in \mathcal{S}$.

To prove the converse, consider
$$
	\mathbf{X} :\overset{d}{=} E \cdot \mathbf1 + \mathbf{S},
$$
where $E$ is an $\text{Exp}(1)$ variable independent of $\mathbf{S}$.  It is trivial to see that \eqref{eq:mpot:std} holds for $\mathbf{X}$.  To show that $\mathbf{Z}$ is a standardized multivariate generalized Pareto distribution, it suffices to show that $\mathbf{X} \in \mathcal{X}$.  Consider the left hand side of \eqref{eq:bm},
\begin{align*}
	\lim_{n\to\infty} \textup{pr}\left\{\max_{i=1,\ldots,n}\mathbf{X}_i - \log(n) \cdot \mathbf1 \le \mathbf{x}\right\}
	=&\, \lim_{n\to\infty} \textup{pr}^n\left\{\mathbf{X} - \log(n) \cdot \mathbf1 \le \mathbf{x}\right\} \\
	=&\, \lim_{n\to\infty} \textup{pr}^n\left\{E \cdot \mathbf1 + \mathbf{S} - \log(n) \cdot \mathbf1 \le \mathbf{x}\right\} \\
	=&\, \lim_{n\to\infty} \left(E_{\mathbf{S}} \left[\textup{pr}_E\left\{E \cdot \mathbf1 + \mathbf{S} - \log(n) \cdot \mathbf1 \le \mathbf{x}\right\} \right] \right)^n \\
	=&\, \lim_{n\to\infty} \left(E_{\mathbf{S}} \left[\textup{pr}_E\left\{E \cdot \mathbf1  \le \mathbf{x} -  \mathbf{S} + \log(n) \cdot \mathbf1\right\} \right] \right)^n \\
	=&\, \lim_{n\to\infty} \left\{E_{\mathbf{S}} \left(\textup{pr}_E\left[E \le \min\left\{\mathbf{x} -  \mathbf{S} + \log(n) \cdot \mathbf1\right\}\right] \right) \right\}^n \\
	=&\, \lim_{n\to\infty} \left[E_{\mathbf{S}} \left\{ 1 - e^{-\log(n) \cdot \mathbf1 + \max(\mathbf{S} - \mathbf{x})} \right\} \right]^n \\
	=&\, \lim_{n\to\infty} \left[E_{\mathbf{S}} \left\{ 1 - \frac{1}{n} \cdot e^{\max(\mathbf{S} - \mathbf{x})} \right\} \right]^n \\
	=&\, \lim_{n\to\infty} \left[1 - \frac{1}{n} \cdot E_{\mathbf{S}} \left\{ e^{\max(\mathbf{S} - \mathbf{x})} \right\} \right]^n \\
	=& \, e^{-E_{\mathbf{S}} \left\{ e^{\max(\mathbf{S} - \mathbf{x})} \right\} } \\
	=&: \, G(\mathbf{x}).
\end{align*}
In particular, $G(x_k) = \exp[-\exp\{-(x_k-\mu_k)\}]$ where $\mu_k=\log\{E(e^{S_k})\}$.  Therefore $\mathbf{X}\in\mathcal{X}$.

\end{proof}


\begin{proof}[Proof of Proposition~\ref{prop:dpot}]

The distribution function of $\mathbf{X} - r \cdot \mathbf1 \left|\{\bar{X} \ge r\} \right.$ can be written as
\begin{align*}
	\textup{pr}(\mathbf{X} - r \cdot \mathbf1 \le \mathbf{z}|\bar{X} \ge r) &= \textup{pr}\{\mathbf{X} - r \cdot \mathbf1 \le \mathbf{z}|\bar{X} \ge r, \max(\mathbf{X})\ge r\} \\
	&= \frac{\textup{pr}\{\mathbf{X} - r \cdot \mathbf1 \le \mathbf{z}, \bar{X} \ge r, \max(\mathbf{X})\ge r\}}{\textup{pr}\{\bar{X} \ge r, \max(\mathbf{X}) \ge r\}} \\
	&= \frac{\textup{pr}\{\mathbf{X} - r \cdot \mathbf1 \le \mathbf{z}, \bar{X} \ge r, \max(\mathbf{X})\ge r\}/\textup{pr}\{\max(\mathbf{X}) \ge r\}}{\textup{pr}\{\bar{X} \ge r, \max(\mathbf{X}) \ge r\}/\textup{pr}\{\max(\mathbf{X}) \ge r\}} \\
	&= \frac{\textup{pr}\{\mathbf{X} - r \cdot \mathbf1 \le \mathbf{z}, \bar{X} \ge r|\max(\mathbf{X})\ge r\}}{\textup{pr}\{\bar{X} \ge r|\max(\mathbf{X}) \ge r\}} \\
	&= \frac{\textup{pr}\{\mathbf{X} - r \cdot \mathbf1 \le \mathbf{z},\bar{X} \ge r|\max(\mathbf{X}) \ge r\}}{\textup{pr}\{(\mathbf{X} - r \cdot \mathbf1)^T\mathbf1\ge 0|\max(\mathbf{X})  \ge r\}}.
\end{align*}
Taking the limit $r\to \infty$ on both sides, 
\begin{align*}
	\lim_{r\to\infty} \textup{pr}(\mathbf{X} - r \cdot \mathbf1 \le \mathbf{z}|\bar{X} \ge r) &= \frac{\lim_{r\to\infty}\textup{pr}\{\mathbf{X} - r \cdot \mathbf1 \le \mathbf{z}, \bar{X} \ge r|\max(\mathbf{X}) \ge r\}}{\lim_{r\to\infty}\textup{pr}\{(\mathbf{X} - r \cdot \mathbf1)^T\mathbf1\ge 0|\max(\mathbf{X}) \ge r\}} \\
	&= \frac{\textup{pr}(\mathbf{Z} \le \mathbf{z}, \mathbf{Z}^T\mathbf1 \ge 0)}{\textup{pr}(\mathbf{Z}^T\mathbf1 \ge 0)} \\
	&= \textup{pr}(\mathbf{Z} \le \mathbf{z}| \mathbf{Z}^T\mathbf1 \ge 0).
\end{align*}
It remains to justify that $\textup{pr}(\mathbf{Z}^T\mathbf1 \ge 0)>0$.  Since the components of $\mathbf{X}$ and hence $\mathbf{Z}$ are asymptotically dependent such that $\textup{pr}(S_j>-\infty)=1$ for $1\le j \le d$, there exists $M>0$ such that $\textup{pr}\{\min(\mathbf{S})>-M\} >0$.  It follows that
\begin{align*}
	\textup{pr}(\mathbf{Z}^T\mathbf1 \ge 0) &= \textup{pr}\{(E\cdot \mathbf1 + \mathbf{S})^T\mathbf1 \ge 0\} \\
	&\ge \textup{pr}\{\min(\mathbf{S})>-M, E>M\} \\
	&= \textup{pr}\{\min(\mathbf{S})>-M\} \cdot \textup{pr}(E>M) >0.
\end{align*}
Therefore $\mathbf{Z^*} \overset{d}= \mathbf{Z} \,|\,{\mathbf{Z}^T\mathbf1 \ge 0}$.

\end{proof}


\begin{proof}[Proof of Proposition~\ref{prop:u}]

Since $\mathbf{Z}$ can be written as $\mathbf{Z} \overset{d}= E \cdot \mathbf1 + \mathbf{T} - \max{(\mathbf{T})} \cdot \mathbf{1}$ and $\mathbf{T} \in \mathbf1^\perp$ such that $\mathbf{T}^T\mathbf1=0$, the conditional event
$$
\{\mathbf{Z}^T \mathbf1 \ge 0\} = \{E \cdot d + \mathbf{T}^T \mathbf1 - \max(\mathbf{T}) \cdot d \ge 0\}  = \{E \cdot d - \max(\mathbf{T}) \cdot d \ge 0\}= \{E -\max(\mathbf{T}) \ge 0\}.
$$
Therefore
$$	
	\mathbf{Z^*} \overset{d}=\left.\{E  - \max(\mathbf{T})\} \cdot \mathbf{1}  + \mathbf{T}\, \right|\, \{E\ge \max(\mathbf{T})\}.
$$
For any $s \ge 0$ and Borel set $B\subseteq \mathbf1^\perp$,
\begin{align*}
	\textup{pr}\{E-\max(\mathbf{T}) \ge s, \mathbf{T} \in B|E \ge \max(\mathbf{T})\}
	&= \frac{\textup{pr}\{E-\max(\mathbf{T}) \ge s, \mathbf{T} \in B\}}{\textup{pr}\{E \ge \max(\mathbf{T})\}} \\
	&= \frac{\int_s^\infty \textup{pr}\{\max(\mathbf{T}) \le t- s, \mathbf{T} \in B\} e^{-t} dt}{\int_0^\infty \textup{pr}\{\max(\mathbf{T}) \le t\} e^{-t} dt}  \\
	&\overset{u=t-s}= \frac{\int_0^\infty \textup{pr}\{\max(\mathbf{T}) \le u, \mathbf{T} \in B\} e^{-(u+s)} du}{\int_0^\infty \textup{pr}\{\max(\mathbf{T}) \le t\} e^{-t} dt} \\
	&= e^{-s}\cdot \frac{\int_0^\infty \textup{pr}\{\max(\mathbf{T}) \le u, \mathbf{T} \in B\} e^{-u} du}{\int_0^\infty \textup{pr}\{\max(\mathbf{T}) \le t\} e^{-t} dt} 
\end{align*}
Take $B= \mathbf{1}^\perp$,
$$
	\textup{pr}\{E-\max(\mathbf{T}) \ge s|E \ge \max(\mathbf{T})\} = e^{-s}.
$$
Take $s=0$, 
$$
	\textup{pr}\{\mathbf{T} \in B|E \ge \max(\mathbf{T})\} = \frac{\int_0^\infty \textup{pr}\{\max(\mathbf{T}) \le u, \mathbf{T} \in B\} e^{-u} dt}{\int_0^\infty \textup{pr}\{\max(\mathbf{T}) \le t\} e^{-t} dt}.
$$
This shows that the conditional distribution of $E-\max(\mathbf{T})\,|\,{E \ge\max(\mathbf{T})}$ is again a unit exponential distribution and $E - \max(\mathbf{T})$ and $\mathbf{T}$ are conditionally independent given $E \ge \max(\mathbf{T})$.
Define
$$
	\mathbf{U} :\overset{d}=\mathbf{T}\,|\,\{\max(\mathbf{T}) \le E\},
$$
then
$$
	\mathbf{Z}^* \overset{d}= E'\cdot\mathbf1 + \mathbf{U}
$$
where $E'$ is an $\text{Exp}(1)$ variable independent of $\mathbf{U}$.

\end{proof}


\begin{proof}[Proof of Corollary~\ref{cor:moments}]

By definition,
\begin{align*}
E\{e^{\max(\mathbf{U})}\} &= E\{e^{\max(\mathbf{T})}\,|\,\max(\mathbf{T}) \le E\} \\
&= \frac{E[e^{\max(\mathbf{T})} \cdot \mathds{1}_{\{\max(\mathbf{T}) \le E\}}]}{E[\mathds{1}_{\{\max(\mathbf{T}) \le E\}}]} \\
&= \frac{\int_0^\infty E[e^{\max(\mathbf{T})} \cdot \mathds{1}_{\{\max(\mathbf{T}) \le s\}}] e^{-s}ds}{\int_0^\infty E[\mathds{1}_{\{\max(\mathbf{T}) \le s\}}] e^{-s}ds} \\
&= \frac{E[e^{\max(\mathbf{T})} \cdot \int_0^\infty \mathds{1}_{\{\max(\mathbf{T}) \le s\}}e^{-s}ds] }{E[\int_0^\infty \mathds{1}_{\{\max(\mathbf{T}) \le s\}}e^{-s}ds] } \\
&= \frac{E\{e^{\max(\mathbf{T})}\cdot e^{-\max(\mathbf{T})}\} }{E\{e^{-\max(\mathbf{T})}\} } \\
&= \frac{1 }{E\{e^{-\max(\mathbf{T})}\} } < \infty.
\end{align*}
The last inequality follows from the fact that $\mathbf{Z}$ has asymptotically dependent components, such that $\textup{pr}(S_k > -\infty) =1$, $k=1,\ldots,d$, and hence $\textup{pr}(T_k < \infty) =1$, $k=1,\ldots,d$.

\end{proof}


\begin{proof}[Proof of Proposition~\ref{prop:t:u}]

Given any Borel set $B \subseteq \mathbf1^\perp$,
\begin{align*}
	\textup{pr}\{\mathbf{U} \in B|\max(\mathbf{U}) = s\}
	 \,=&   \,\textup{pr}\{\mathbf{T} \in B| \max(\mathbf{T}) = s, \max(\mathbf{T}) \le E\}\\
	 \,=&   \,\textup{pr}\{\mathbf{T} \in B| \max(\mathbf{T}) = s, E\ge s\}\\
	 \,=&   \,\textup{pr}\{\mathbf{T} \in B| \max(\mathbf{T}) = s\}.
\end{align*}
Hence
$$
	\mathbf{U}\,|\,\{\max(\mathbf{U})=s\} \overset{d} = \mathbf{T}\,|\,\{\max(\mathbf{T})=s\}.
$$
For any $s \ge 0$, 
\begin{align*}
\textup{pr}\{\max(\mathbf{U}) \le s\}  \,=&  \,\textup{pr}\{\max(\mathbf{T}) \le s|\max(\mathbf{T}) \le E\} \\
=& \, \frac{\textup{pr}\{\max(\mathbf{T}) \le s,\max(\mathbf{T}) \le E\}}{ \textup{pr}\{\max(\mathbf{T}) \le E\} } \\
=& \, \frac{\int_0^s \textup{pr}\{\max(\mathbf{T}) \le t\} e^{-t} dt + \int_s^\infty \textup{pr}\{\max(\mathbf{T}) \le s\} e^{-t} dt}{\int_0^\infty \textup{pr}\{\max(\mathbf{T}) \le t\} e^{-t} dt} \\
\overset{u=e^{-t}}=&  \,\frac{\int_0^s \textup{pr}\{\max(\mathbf{T}) \le t\} e^{-t} dt + \textup{pr}\{\max(\mathbf{T}) \le s\} \int_s^\infty e^{-t} dt}{\int_0^1 \textup{pr}\{\max(\mathbf{T}) \le -\log(u)\} du} \\
=& \, \frac{\int_0^s \textup{pr}\{\max(\mathbf{T}) \le t\} e^{-t} dt + e^{-s} \cdot \textup{pr}\{\max(\mathbf{T}) \le s\} }{\int_0^1 \textup{pr}\left\{e^{-\max(\mathbf{T})} \ge u \right\} du} \\
=& \, \frac{\int_0^s \textup{pr}\{\max(\mathbf{T}) \le t\} e^{-t} dt + e^{-s} \cdot \textup{pr}\{\max(\mathbf{T}) \le s\} }{E\left\{e^{-\max(\mathbf{T})} \right\}}.
\end{align*}
This proves Equation \eqref{eq:maxu} of the proposition.

To prove Equation \eqref{eq:maxt}, it suffices to show that the random variable $\max(\mathbf{T})$ defined by \eqref{eq:maxt} satisfies \eqref{eq:maxu}.  Let $\max(\mathbf{T})$ be the random variable defined by distribution function,
$$
\textup{pr}\{\max(\mathbf{T}) \le s\} = \frac{e^{s} \textup{pr}\{\max(\mathbf{U}) \le s\} - \int_0^s \textup{pr}\{\max(\mathbf{U}) \le  t \}  e^{t} dt}{E\left\{e^{\max(\mathbf{U})}\right\}}, \quad s \ge 0.
$$
By definition
$$
	\textup{pr}\{\max(\mathbf{T}) \le s\} \propto e^{s} \textup{pr}\{\max(\mathbf{U}) \le s\}-\int_0^s \textup{pr}\{\max(\mathbf{U}) \le  t \}  e^{t} dt.
$$
Plug it in the numerator of \eqref{eq:maxu}, which is,
$$
	\int_0^s \textup{pr}\{\max(\mathbf{T}) \le t\}  e^{-t} dt + e^{-s} \textup{pr}\{\max(\mathbf{T}) \le s\},
$$
we obtain
\begin{align*}
	& \int_0^s \textup{pr}\{\max(\mathbf{T}) \le t\}  e^{-t} dt + e^{-s} \textup{pr}\{\max(\mathbf{T}) \le s\}\\
	=&\, e^{-s} \left[ e^{s} \textup{pr}\{\max(\mathbf{U}) \le s\}-\int_0^s \textup{pr}\{\max(\mathbf{U}) \le  t \}  e^{t} dt\right] \\
	& \quad +\int_0^s \left[ e^{t} \textup{pr}\{\max(\mathbf{U}) \le t\} - \int_0^t \textup{pr}\{\max(\mathbf{U}) \le  u \}  e^{u} du\right]  e^{-t} dt \\
	 =&\,  \underbrace{\textup{pr}\{\max(\mathbf{U}) \le s\}}_{\text{Term I}} - \underbrace{e^{-s}\int_0^s \textup{pr}\{\max(\mathbf{U}) \le  t \}  e^{t} dt}_{\text{Term II}}  +  \underbrace{\int_0^s \textup{pr}\{\max(\mathbf{U}) \le t\}dt}_{\text{Term III}} \\
	 &\quad - \underbrace{\int_0^s \int_0^t \textup{pr}\{\max(\mathbf{U}) \le  u \}  e^{u-t} dudt}_{\text{Term IV}}.
\end{align*}
Consider Term IV, 
\begin{align*}
	\text{Term IV} =& \int_0^s \int_0^t \textup{pr}\{\max(\mathbf{U}) \le  u  \}  e^{u-t} dudt \\
	=& \int_0^s \int_u^s \textup{pr}\{\max(\mathbf{U}) \le  u  \}  e^{u-t} dtdu\\
	=& \int_0^s  \textup{pr}\{\max(\mathbf{U}) \le  u  \}  e^{u} (e^{-u} - e^{-s}) du\\
	=& \int_0^s  \textup{pr}\{\max(\mathbf{U}) \le  u \}du - e^{-s}\int_0^s  \textup{pr}\{\max(\mathbf{U}) \le  u  \}  e^{u} du\\
	=& \,\text{Term III} - \text{Term II}.
\end{align*}
Hence 
$$
	\int_0^s \textup{pr}\{\max(\mathbf{T}) \le t\}  e^{-t} dt + e^{-s} \textup{pr}\{\max(\mathbf{T}) \le s\}
 \propto \text{Term I} =\textup{pr}\{\max(\mathbf{U}) \le  s\},
$$
which gives \eqref{eq:maxu}, hence proving \eqref{eq:maxt}.

\end{proof}


\begin{proof}[Proof of Corollary~\ref{cor:t:u:2}]

The result follows from Proposition 4 by taking the derivatives of the distribution functions with respect to $s$.

\end{proof}


\begin{proof}[Proof of Corollary~\ref{cor:t:u:1}]

For any $x>0$ and $r>x$,
\begin{align*}
& \, \textup{pr}\left\{\max(\mathbf{U})\le x|{\max(\mathbf{U}) \ge r-E} \right\} \\
=&\,  \frac{\textup{pr}\left\{r-E < \max(\mathbf{U})\le x \right\}}{\textup{pr}\left\{\max(\mathbf{U})\ge r - E\right\}} \\
	=&\, \frac{\int_{r-x}^\infty \textup{pr}\left\{r-t < \max(\mathbf{U})\le x \right\} e^{-t}dt}{\int_{0}^\infty \textup{pr}\left\{\max(\mathbf{U})\ge r - t\right\}e^{-t} dt} \\
	=&\, \frac{\int_{r-x}^\infty \textup{pr}\left\{\max(\mathbf{U})\le x \right\} e^{-t}dt - \int_{r-x}^\infty \textup{pr}\left\{\max(\mathbf{U})\le r-t \right\} e^{-t}dt}{\int_{0}^\infty \textup{pr}\left\{\max(\mathbf{U})\ge r - t\right\}e^{-t} dt} \\
	\overset{u=r-t}=&\frac{\textup{pr}\left\{\max(\mathbf{U})\le x \right\} e^{-r+x} - \int_{-\infty}^x \textup{pr}\left\{\max(\mathbf{U})\le u \right\}e^{u-r}du}{\int_{-\infty}^r \textup{pr}\left\{\max(\mathbf{U})\ge u\right\}e^{u-r} du} \\
	=&\, \frac{\textup{pr}\left\{\max(\mathbf{U})\le x \right\} e^{x} - \int_{0}^x \textup{pr}\left\{\max(\mathbf{U})\le u \right\} e^{u}du}{\int_{-\infty}^r \textup{pr}\left\{\max(\mathbf{U})\ge u\right\}e^{u } du } \\
	\overset{r\to\infty} \to&\frac{\textup{pr}\left\{\max(\mathbf{U})\le x \right\} e^{x} - \int_{0}^x \textup{pr}\left\{\max(\mathbf{U})\le u \right\} e^{u}du}{\int_{-\infty}^\infty \textup{pr}\left\{\max(\mathbf{U})\ge u\right\}e^{u } du } \\
	 = &\, \frac{\textup{pr}\left\{\max(\mathbf{U})\le x \right\} e^{x} - \int_{0}^x \textup{pr}\left\{\max(\mathbf{U})\le u \right\} e^{u}du}{\int_{-0}^\infty \textup{pr}\left\{e^{\max(\mathbf{U})}\ge e^u\right\}d e^u } \\
	= &\, \frac{\textup{pr}\left\{\max(\mathbf{U})\le x \right\} e^{x} - \int_{0}^x \textup{pr}\left\{\max(\mathbf{U})\le u \right\} e^{u}du}{E\left\{e^{\max(\mathbf{U})}\right\} } \\
	=&\,  \textup{pr}\left\{\max(\mathbf{T})\le x \right\}.
\end{align*}

\end{proof}


\begin{proof}[Proof of Proposition~\ref{prop:rv}]

The proof follows the same lines as that of Proposition~\ref{prop:rep:mgp}, replacing $\mathbf{S}$ with $\mathbf{U}$.

\end{proof}


\begin{proof}[Proof of Proposition~\ref{prop:hr}]

Let 
$$
	\Sigma := -\frac{1}{2}\left(I - d^{-1}{\mathbf{1}\mathbf{1}^T}\right) \Gamma \left(I -d^{-1}{\mathbf{1}\mathbf{1}^T}\right)
$$
be the covariance matrix corresponding to the variogram $\Gamma$ that satisfies $\Sigma \mathbf1=\mathbf0$.  In the case where $\Sigma$ is of rank $(d-1)$, \cite{HES2024} derived the density of the standardized multivariate generalized Pareto distribution $\mathbf{Z}$, which gave this result.

Assume that $\Sigma$ is of rank lower than $(d-1)$.  Define
$$
	\Gamma_m:= \Gamma + \frac2m (\mathbf1\mathbf1^T - I), \quad m =1,2,\ldots
$$
Then each $\Gamma_m$ is the variogram corresponding to the covariance matrix $\Sigma_m:=\Sigma + m^{-1} I$, which is of full rank, and $\Gamma_m \to \Gamma$ as $m\to\infty$.

Let $G_m, \mathbf{Z}_m, \mathbf{Z}_m^*,\mathbf{U}_m$ be the generalized multivariate extreme value distribution, standardized multivariate generalized Pareto distirbution, diagonal multivariate generalized Pareto distribution and profile random vector of the H\"usler-Reiss model parametrized by $\Gamma_m$, respectively.  Let $G, \mathbf{Z}, \mathbf{Z}^*,\mathbf{ U}$ be that of the H\"usler-Reiss model parametrized by $\Gamma$.  
The generalized multivariate extreme value distribution $G$ can be expressed as a function of the variogram $\Gamma$, see \cite{HR1989}, such that the convergence of variogram $\Gamma_m \to \Gamma$ leads to the convergence in in distribution of the generalized multivariate extreme value distributions,
$$
	G_m(\mathbf{x}) \to G(\mathbf{x}), \quad \mathbf{x} \in \mathbb{R}^d.
$$
From the expression of the distribution function of the standardized multivariate generalized Pareto distribution $\mathbf{Z}$ in \eqref{eq:mgpd:df}, this implies that $\mathbf{Z}_m \overset{d}\to \mathbf{Z}$ and hence $ \mathbf{Z}_m^* \overset{d}\to \mathbf{Z}^*$. By definition, the profile random vector $\mathbf{U}$ is the projection of $\mathbf{Z}$ onto $\mathbf1^\perp$ and can therefore be obtained from $\mathbf{Z}$ via the continuous mapping
$$
	\mathbf{U} \overset{d}= \mathbf{Z}^* - \bar{Z}^* \cdot \mathbf1.
$$
Consequently,
$$
	 \mathbf{U}_m\overset{d}\to \mathbf{U}.
$$
Define $\boldsymbol\mu:= -\left\{\text{diag}(\Sigma) - \overline{\text{diag}(\Sigma)} \cdot \mathbf1\right\}/2$ and $\boldsymbol\mu_m:= -\left\{\text{diag}(\Sigma_m) - \overline{\text{diag}(\Sigma_m)} \cdot \mathbf1\right\}/2$.
Since
$$
	\mathbf{U}_m = N(\boldsymbol\mu_m,\Sigma_m) \overset{d}\to N(\boldsymbol\mu,\Sigma),
$$
we obtain
$$
	\mathbf{U} \sim N(\boldsymbol\mu,\Sigma).
$$

\end{proof}


\bibliographystyle{chicago}
\bibliography{ref}

\end{document}